\documentclass[letterpaper, 10pt, conference]{ieeeconf}
\IEEEoverridecommandlockouts 

\usepackage{amsmath,amssymb}
\usepackage{epsfig,curves,multirow,subfigure,url,array,warmread}
\usepackage[all,import]{xy}

\newcommand{\norm}[1]{\ensuremath{\left\| #1 \right\|}}
\newcommand{\bracket}[1]{\ensuremath{\left[ #1 \right]}}
\newcommand{\braces}[1]{\ensuremath{\left\{ #1 \right\}}}
\newcommand{\parenth}[1]{\ensuremath{\left( #1 \right)}}
\newcommand{\refeqn}[1]{(\ref{eqn:#1})}
\newcommand{\reffig}[1]{Fig. \ref{fig:#1}}

\newcommand{\deriv}[2]{\ensuremath{\frac{\partial #1}{\partial #2}}}
\newcommand{\T}{\ensuremath{\mathrm{T}}}
\newcommand{\SO}{\ensuremath{\mathrm{SO(3)}}}
\newcommand{\so}{\ensuremath{\mathfrak{so}(3)}}
\newcommand{\SE}{\ensuremath{\mathrm{SE(3)}}}
\newcommand{\se}{\ensuremath{\mathfrak{se}(3)}}

\renewcommand{\Re}{\ensuremath{\mathbb{R}}}
\renewcommand{\S}{\ensuremath{\mathbb{S}}}

\title{\LARGE \bf
A Combinatorial Optimal Control Problem for Spacecraft Formation Reconfiguration}

\author{ \parbox{3 in}{\centering Taeyoung Lee\authorrefmark{1}\authorrefmark{2}, N. Harris McClamroch\authorrefmark{2}\\
         Department of Aerospace Engineering\\
         University of Michigan, Ann Arbor, MI 48109\\
         {\tt\small \{tylee, nhm\}@umich.edu}}
         \hspace*{ 0.5 in}
         \parbox{3 in}{\centering Melvin Leok\authorrefmark{1}\\
         Department of Mathematics\\
        Purdue University, West Lafayette, IN 47907\\
         {\tt\small mleok@math.purdue.edu}}
        \thanks{\textsuperscript{\footnotesize\ensuremath{*}}This research has been supported in part by NSF under grant DMS-0504747, and by a grant from the Rackham Graduate School, University of Michigan.}
        \thanks{\textsuperscript{\footnotesize\ensuremath{\dagger}}This research has been supported in part by NSF under grant ECS-0244977.}
}

\begin{document}
\allowdisplaybreaks
\maketitle \thispagestyle{empty} \pagestyle{empty}

\begin{abstract}
We consider a spacecraft formation reconfiguration problem in the case of identical spacecraft. This introduces in the optimal reconfiguration problem a permutational degree of freedom, in addition to the choice of individual spacecraft trajectories. We approach this using a coupled combinatorial and continuous optimization framework, in which the inner loop consists of computing the costs associated with a particular assignment by using a geometrically exact and numerically efficient discrete optimal control method based on Lie group variational integrators. In the outer optimization loop, combinatorial techniques are used to determine the optimal assignments based on the costs computed in the inner loop. The proposed method is demonstrated on the optimal reconfiguration problem for 5 identical spacecraft to go from an inline configuration to one equally spaced on a circle.
\end{abstract}

\section{Introduction}

The optimal control of spacecraft formations has received increased interest due to the NASA Terrestrial Planet Finder Project and the ESA Darwin Project. The objective is to use multiple spacecraft for cooperative missions such as long base-line interferometers, magnetosphere studies, and space-based communication networks. To accomplish various goals efficiently, it is often required to reconfigurate a formation during a mission. Since each spacecraft has limited fuel, it is important that the formation reconfiguration maneuvers are achieved with minimum fuel expenditure.

Formation reconfiguration can be classified into two types: each spacecraft is required to be transferred into a specified location in the desired reconfigured formation, while in the other case, a specified location in the desired formation can be occupied by any spacecraft of a particular type~\cite{WanHad.JAS99}. In general, a formation is composed of identical spacecraft or groups of spacecraft at the same type, and the total fuel consumption depends on the permutations of the formation reconfigurations as well as the maneuver of each spacecraft to a specified location.

In this paper, we study an optimal spacecraft formation problem integrated with an integer/combinatorial optimization approach for the assignment. Usually in combinatorial optimization problems for multiple agents, the dynamics of each agent is either ignored or simplified into an analytic model such as a kinematics equation or a double integrator~\cite{SavBulFra.CDC06}. Here, we assume that each spacecraft evolves on the special Euclidean group $\SE$, including both translational dynamics and rotational attitude dynamics under a central gravitational potential. Thus, finding the optimal control forces and moments on a spacecraft assigned to a fixed desired location is demanding even if the combinatorial assignment optimization problem is not considered. This is an interesting and challenging problem since it requires combining an integer/combinatorial optimization approach and an optimal control method over the non-trivial dynamics of spacecraft on $\SE$.

There has been some work in the literature on combinatorial optimization for spacecraft formation. The costs for all possible assignments are directly compared for spacecraft moving along a straight path in~\cite{WanHad.JAS99}. This requires a large computational effort since optimal control problems associated with $n!$ assignments have to be solved for a formation of $n$ spacecraft. In~\cite{GuiSch.JGCD06}, the special structure of a Hamiltonian system is utilized to expedite finding the solutions of optimization problems with varying boundary conditions. But, this requires a solution of the Hamilton-Jacobi partial differential equation. A stochastic optimization technique is used in~\cite{HosAtk.AIAA05}.

To approximate the cost matrix used in the combinatorial assignment problem, we use the cost entries which have been explicitly computed, and the sensitivities of the cost to construct approximations to the remaining entries. The solution of the optimal control problem for each spacecraft is based on a structure-preserving numerical integrator referred to as a Lie group variational integrator~\cite{CMAME05}. Combined with an indirect optimization method, the Lie group variational integrator provides a geometrically exact but numerically efficient numerical optimization method for the dynamics of a rigid body~\cite{CDC06.opt}. The combinatorial optimization scheme for spacecraft formations that we present has the following important features: (1) dynamics of each spacecraft is nontrivial, (2) a discrete combinatorial optimization on a permutation group is explicitly integrated into the continuous optimal control problems, and (3) the problem is formulated and solved in a discrete time space using a Lie group variational integrator for overall computational accuracy and efficiency.

This paper is organized as follows. Computational approaches to solve an optimal control problem for a single spacecraft are summarized in Section \ref{sec:lgvi} and \ref{sec:opt}. Based on these results, a combinatorial optimization approach is developed in Section \ref{sec:comb}, which is followed by a numerical example in Section \ref{sec:ne}.

\section{Lie Group Variational Integrator}\label{sec:lgvi}
The configuration space for the translational and rotational motion of a rigid body is the special Euclidean group, $\SE=\Re^3\,\textcircled{s}\,\SO$. We identify the cotangent bundle $\T^*\SE$ with $\SE\times\se^*$ by left translation, and we identify $\se^*$ with $\Re^6$ by an isomorphism between $\Re^6$ and $\se$, and the isomorphism between $\se$ and $\se^*$ induced by the standard inner product on $\Re^6$. We denote the attitude, position, angular momentum, and linear momentum of the rigid body by $(R,x,\Pi,\gamma)\in\T^*\SE$.

The continuous equations of motion are given by
\begin{gather}
\dot{x}=\frac{\gamma}{m},\label{eqn:xdot}\\
\dot{\gamma}=f+u^f,\label{eqn:gamdot}\\
\dot{R}=RS(\Omega),\label{eqn:Rdot}\\
\dot{\Pi}+\Omega\times\Pi = M+u^m,\label{eqn:Pidot}
\end{gather}
where $\Omega\in\Re^3$ is the angular velocity, and $u^f,u^m\in\Re^3$ are the control force in the inertial frame and the control moment in the body fixed frame, respectively. The constant mass of the rigid body is $m\in\Re$, and $J\in\Re^{3\times 3}$ denotes the moment of inertia, i.e. $\Pi=J\Omega$. The map $S(\cdot):\Re^3\mapsto\so$ is an isomorphism between $\so$ and $\Re^3$ defined by the condition $S(x)y=x\times y$ for all $x,y\in\Re^3$.

We assume that the potential is dependent on the position and the attitude; $U(\cdot):\SE\mapsto\Re$. The corresponding force and the moment due to the potential are given by
\begin{align}
f & =-\deriv{U}{x},\label{eqn:f}\\
M & =r_{1}\times u_{r_1}+r_{2}\times u_{r_2}+r_{3}\times u_{r3},\label{eqn:M}
\end{align}
where $r_i,u_{r_i}\in\Re^3$ are the $i$th row vector of $R$ and $\deriv{U}{R}$, respectively.

Since the dynamics of a rigid body has the structure of a Lagrangian or Hamiltonian system, they are symplectic, momentum and energy preserving. These geometric features determine the qualitative behavior of the rigid body dynamics, and they can serve as a basis for theoretical study of rigid body dynamics.

In contrast, the most common numerical integration methods, including the widely used (non-symplectic) explicit Runge-Kutta schemes, preserve neither the Lie group structure nor these geometric properties. In addition, standard Runge-Kutta methods fail to capture the energy dissipation of a controlled system accurately~\cite{MarWes.AN01}. Additionally, if we integrate \refeqn{Rdot} by a typical Runge-Kutta scheme, the quantity $R^T R$ inevitably drifts from the identity matrix as the simulation time increases. It is often proposed to parameterize rotations by Euler angles or unit quaternions. However, Euler angles are not global expressions of the attitude since they have associated singularities. Unit quaternions do not exhibit singularities, but are constrained to lie on the unit three-sphere $\S^3$, and general numerical integration methods do not preserve the unit length constraint. Therefore, quaternions have the same numerical drift problem. Renormalizing the quaternion vector at each step tends to break the conservation properties. Furthermore, unit quaternions, which are diffeomorphic to $\mathrm{SU(2)}$, double cover $\SO$. So there are inevitable ambiguities in expressing the attitude.

In~\cite{CMAME05}, Lie group variational integrators are constructed by explicitly adapting  Lie group methods~\cite{IserMun.AN00} to the discrete variational principle~\cite{MarWes.AN01}. They have the desirable property that they are symplectic and momentum preserving, and they exhibit good energy behavior for an exponentially long time period. They also preserve the Euclidian Lie group structure without the use of local charts, reprojection, or constraints. These geometrically exact numerical integration methods yield highly efficient and accurate computational algorithms for rigid body dynamics, and avoid singularities and ambiguities.

Using the results presented in~\cite{CMAME05}, a Lie group variational integrator on $\SE$ for equations \refeqn{xdot}--\refeqn{Pidot} is given by
\begin{gather}
x_{k+1} = x_k +\frac{h}{m} \gamma_k + \frac{h^2}{2m}\parenth{f_k+u^f_k},\label{eqn:xkp}\\
\gamma_{k+1} = \gamma_k +\frac{h}{2}\parenth{f_k+u^f_k}+\frac{h}{2}\parenth{f_{k+1}+u^f_{k+1}},\label{eqn:gamkp}\\
h S(\Pi_k + \frac{h}{2}\parenth{M_k+u^m_k}) = F_k J_d -J_d F_k^T,\label{eqn:findf}\\
R_{k+1} = R_k F_k,\label{eqn:Rkp}\\
\Pi_{k+1} = F_k^T\Pi_k +\frac{h}{2}F_k^T \parenth{M_k+u^m_k} +
\frac{h}{2}\parenth{M_{k+1}+u^m_{k+1}},\label{eqn:Pikp}
\end{gather}
where the subscript $k$ denotes the $k$-th step for a fixed integration step size $h\in\Re$. The matrix $J_d\in\Re^{3\times 3}$ is a nonstandard moment of inertia matrix defined by $J_d=\frac{1}{2}\mathrm{tr}\!\bracket{J}I_{3\times 3}-J$. The matrix $F_k\in\SO$ denotes the relative attitude between adjacent integration steps.

For given $(R_k,x_k,\Pi_k,\gamma_k)$ and control input, \refeqn{findf} is solved to find $F_k$. Then $(R_{k+1},x_{k+1})$ are obtained by \refeqn{Rkp} and \refeqn{xkp}. Using \refeqn{f} and \refeqn{M}, $(f_{k+1},M_{k+1})$ are computed, and they are used to find $(\Pi_{k+1},\gamma_{k+1})$ by \refeqn{Pikp} and\refeqn{gamkp}. This yields a map $(R_k,x_k,\Pi_k,\gamma_k)\mapsto(R_{k+1},x_{k+1},\Pi_{k+1},\gamma_{k+1})$, and this process is repeated. The only implicit part is \refeqn{findf}. The actual computation of $F_k$ is done in the Lie algebra $\so$ of dimension 3, and the rotation matrices are updated by multiplication. As this approach does not involve the component-wise integration of the kinematics equation \refeqn{Rdot}, there is no excessive computational burden. It can be shown that this integrator has second-order accuracy.

One of the distinct features of the Lie group variational integrator is that it preserves both the symplectic property and the Lie group structure of the rigid body dynamics. As such, it exhibits substantially improved computational accuracy and efficiency compared with other geometric integrators that preserve only one of these properties, that is symplectic Runge-Kutta methods that do not preserve Lie group structure or non-symplectic Lie group methods~\cite{CMDA06}. The symplectic property of numerical integrators is important even in the case of controlled dynamics, since the dissipation rate of the total energy is typically computed inaccurately by non-symplectic integrators~\cite{MarWes.AN01}.

\section{Optimal Control of a Rigid Body on \SE}\label{sec:opt}
We first summarize a computational approach to solve the optimal control problem for a single rigid body in which the translational dynamics and the rotational attitude dynamics are coupled~\cite{CDC06.opt}. This approach is extended to solve a combinatorial optimal formation problem for multiple rigid bodies in Section \ref{sec:comb}.

\subsection{Problem Formulation}
An optimal control problem is formulated for the maneuver of a rigid body from a given initial configuration $(R_0,x_0,\Pi_0,\gamma_0)$ to a desired configuration $(R_N^d,x_N^d,\Pi_N^d,\gamma_N^d)$ during a given maneuver time $N$. Control inputs are parameterized by their value at each time step. The performance index is the square of the weighted $l_2$ norm of the control inputs.
\begin{gather*}
\text{given: } (x_0,\gamma_0,R_0,\Pi_0),\,(x_d,\gamma_d,R_d,\Pi_d),\,N,\\
\min_{u} \mathcal{J}=\sum_{k=0}^{N-1}
\braces{\frac{h}{2}(u^f_{k+1})^TW_fu^f_{k+1}+
\frac{h}{2}(u^m_{k+1})^TW_mu^m_{k+1}},\\
\text{such that } (x_N,\gamma_N,R_N,\Pi_N)=(x_d,\gamma_d,R_d,\Pi_d),\\
\text{subject to discrete equations of motion
\refeqn{xkp}--\refeqn{Pikp},}
\end{gather*}
where $W_f,W_m\in\Re^{3\times 3}$ are symmetric positive definite
matrices.

\subsection{Computational Approach}
We solve this optimal control problem using an indirect method; necessary conditions for optimality are obtained by using variational expressions that respect the geometry of the configuration space, and the corresponding two point boundary value problem is solved by using the shooting method. Here we use a modified version of the discrete equations
of motion with first order accuracy, because it yields a compact form for the necessary conditions.

Define an augmented performance index as
\begin{align*}
\mathcal{J}_a =
\sum_{k=0}^{N-1}&\frac{h}{2}(u^f_{k+1})^TW^fu^f_{k+1}+
\frac{h}{2}(u^m_{k+1})^TW^mu^m_{k+1}\nonumber\\
& +\lambda_k^{1,T}\braces{-x_{k+1}+x_k+\frac{h}{m}\gamma_k}\nonumber\\
& +\lambda_k^{2,T}\braces{-\gamma_{k+1} + \gamma_k+hf_{k+1}+hu^f_{k+1}}\nonumber\\
& +\lambda_k^{3,T}S^{-1}\!\parenth{\mathrm{logm}(F_k-R_{k}^TR_{k+1})}\nonumber\\
& +\lambda_k^{4,T}\braces{-\Pi_{k+1} + F_k^T \Pi_k +
h\parenth{M_{k+1}+u_{k+1}^m}},
\end{align*}
where $\lambda_k^{i}\in \Re^3$ are Lagrange multipliers. The constraint \refeqn{findf} is considered implicitly using a constrained variation. Setting $\delta\mathcal{J}_a=0$ for all variations,  we obtain necessary conditions for optimality as follows.
\begin{gather}
u^f_{k+1} = -W_{f}^{-1}\lambda_{k}^2,\label{eqn:ufkp}\\
u^m_{k+1} = -W_{m}^{-1}\lambda_{k}^4,\label{eqn:umkp}\\
\lambda_{k} = A_{k+1}^T \lambda_{k+1},\label{eqn:updatelam}
\end{gather}
$\lambda_k=[\lambda_k^1;\lambda_k^2;\lambda_k^3;\lambda_k^4]\in\Re^{12}$,
and $A_k\in\Re^{12\times 12}$ is suitably defined in terms of $(R_k,x_k,\Pi_k,\gamma_k)$. Together with the discrete equations of motion, this
yields a map $\braces{(R_k,x_0,\Pi_k,\gamma_k),\lambda_k}\mapsto
\braces{(R_{k+1},x_{k+1},\Pi_{k+1},\gamma_{k+1}),\lambda_{k+1}}$.

The necessary conditions for optimality are expressed in terms of a two point boundary problem on $\T^{*}\SE$ and its dual. This problem is to find the optimal discrete flow, multiplier, and control inputs to satisfy the equations of motion, optimality conditions \refeqn{ufkp},\refeqn{umkp}, multiplier equations \refeqn{updatelam}, and boundary conditions simultaneously.

We use the shooting method~\cite{Bry.BK75}. A nominal solution satisfying all of the necessary conditions except the boundary conditions is chosen. The unspecified initial multiplier is updated by successive linearization so as to satisfy the specified terminal boundary conditions in the limit. The optimality conditions \refeqn{ufkp} and \refeqn{umkp} are substituted into the equations of motion and the multiplier equations. The sensitivities of the specified terminal boundary conditions with respect to the unspecified initial multiplier conditions is obtained by a linear analysis.

Let $z_k\in\Re^{12}$ be the variation of the state given by $z_k=[\xi_k;\delta x_k;\delta\Pi_k;\delta\gamma_k]$, where $\zeta_k\in\Re^3$ denotes the variation of the rotation matrix as $\delta R_k=\frac{d}{d\epsilon} \big|_{\epsilon=0}R_k\exp S(\zeta_k)=R_k S(\zeta_k)$, and variations for other variables are defined in the usual sense. The linearized equations of motion and the linearized multiplier equation can be written as
\begin{align}
z_{k+1} & = A_k z_k + \mathcal{A}^{12} \delta\lambda_k,\label{eqn:zkp}\\
\delta \lambda_k & = \mathcal{A}_{k+1}^{21} z_{k+1} + A_{k+1}^T \delta
\lambda_{k+1},\label{eqn:dellamkp}
\end{align}
where $\mathcal{A}^{12}_k,\mathcal{A}_{k+1}^{21}\in\Re^{12\times 12}$ can be computed explicitly. The solution of the linear equations \refeqn{zkp} and \refeqn{dellamkp} can be obtained as
\begin{align}
\begin{bmatrix}z_{k}\\\delta\lambda_{k}\end{bmatrix}
= \begin{bmatrix} \Phi_k^{11} & \Phi_k^{12}\\\Phi_k^{21} &
\Phi_k^{22}\end{bmatrix}
\begin{bmatrix}z_{0}\\\delta\lambda_{0}\end{bmatrix},\label{eqn:varopt}
\end{align}
where $\Phi_k^{ij}\in\Re^{12\times 12}$.

For the given two point boundary value problem, $z_0=0$ since the initial condition is fixed, and $\lambda_N$ is free.  Thus,
\begin{align}
z_N = \Phi^{12}_N \delta\lambda_0.
\end{align}
The matrix $\Phi^{12}_N$ represents the sensitivity of the specified terminal boundary conditions with respect to the unspecified initial multipliers. Using this sensitivity, an initial guess of the unspecified initial conditions is iterated to satisfy the specified terminal conditions in the limit. Any type of Newton iteration can be applied. We use a line search with backtracking algorithm, referred to as Newton-Armijo iteration in~\cite{Kel.BK95}: the outer loop computes the sensitivity derivatives to obtain the Newton search direction, and the inner loop performs a line search to find the largest step size along the given search direction.

\subsection{Properties of Computational Approach}

The key feature of this computational approach for the optimal control problem of a single rigid body is that it is discretized from the problem definition level using the Lie group variational integrator. This is in contrast to obtaining continuous time necessary conditions, which are discretized to numerically solve the two point boundary value problem. In this computational approach for the optimal control problem, the  discrete necessary conditions for optimality are obtained by a variational principle.

The main advantage of the shooting method is that the number of iteration variables, the initial Lagrange multipliers, is small. In other approaches, an initial guess of a control input history or multiplier history are iterated, so the number of optimization parameters is proportional to the number of discrete time steps. The difficulty is that the extremal solutions are sensitive to small changes in the unspecified initial multiplier values. The nonlinearities also make it hard to construct an accurate estimate of sensitivity, perhaps resulting in numerical ill-conditioning.

Here, the discrete necessary conditions for optimality preserve the geometric structure of the optimal control problem. Thus, there is no geometrical error introduced by the numerical integration algorithm itself. It turns out that, combined with the shooting method, this computational approach provides a geometrically exact and numerically efficient solution to this highly nonlinear, non-convex rigid body optimal control problem~\cite{CDC06.opt,ACC07.opt}. This is used as a basic tool for the combinatorial optimization problem for multiple rigid bodies.

\section{Optimal Formation Control of Rigid Bodies}\label{sec:comb}

\subsection{Problem Formulation}

We study an optimal formation control problem of $n$ identical rigid bodies where the maneuver of each body is described by \refeqn{xkp}--\refeqn{Pikp}. The objective is to find the optimal control forces and moments for each rigid body such that the group moves from a given initial configuration $(R_0^i,x_0^i,\Pi_0^i,\gamma_0^i)$ for $i\in\braces{1,2,\ldots,n}$ to a desired target $\mathcal{T}\in\Re^{3n}$ during a given maneuver time $N$, where the superscript $i$ denotes the $i$-th rigid body.

More precisely, we assume that the $n$ desired positions $\braces{x_d^i(\theta)}_{i=1}^n$, at which all rigid bodies are located at the terminal maneuver time, are given as functions of parameters $\theta\in\Re^l$. The desired attitude, the linear momentum, and the angular momentum at the terminal time, $(R_d,\Pi_d,\gamma_d)$, are assumed to be fixed and to be the same for all rigid bodies.

Since all rigid bodies are identical, there are $n!$ possible combinatorial assignments for $n$ rigid bodies to these $n$ desired locations. Let $\braces{a_{ij}}$ be a $n\times n$ matrix composed of binary elements $\braces{0,1}$, referred to as an assignment or a permutation matrix. Each element of the assignment matrix $a_{ij}$ represents the assignment of the $i$-th rigid body to the $j$-th desired terminal position $x_d^j$. If $a_{ij} = 1$, the $i$-th rigid body is assigned to the $j$-th node, and if $a_{ij} = 0$, the $i$-th rigid body is not assigned to the $j$-th node. Thus, the assignment is valid if $\sum_{j=1}^na_{ij}=\sum_{i=1}^n a_{ij}=1$. The assignment matrix can be equivalently expressed as a set $A=\braces{(i, j) | a_{ij} = 1}$, and the particular desired points assigned by the $i$-th rigid body for an assignment $A$ is denoted by $A_i\in\braces{1,\ldots,n}$. In other words, for an assignment $A$, the $i$-th rigid body is assigned to the $A_i$-th desired location, $x^{A_i}_d$.

The target is defined in terms of a parameter $\theta$ and an assignment $A$ as follows.
\begin{align*}
    \mathcal{T}(\theta,A)=\braces{x^{A_i}_d(\theta)}_{i=1}^n \in \Re^{3n}.
\end{align*}
Thus, for a given parameter $\theta\in\Re^l$ and a given assignment $A\in S_n$, the terminal boundary conditions for all rigid bodies are completely determined.

The performance index is the sum of the squares of the weighted $l_2$ norms of the control inputs. The optimal control problem for a formation of $n$ rigid bodies is formulated as
\begin{gather*}
\text{given: } \braces{(x_0^i,\gamma_0^i,R_0^i,\Pi_0^i)}_{i=1}^n,\,\braces{x_d^i(\theta)}_{i=1}^n,R_d,\Pi_d,\gamma_d),\,N,\\
\min_{u,\theta,A} \mathcal{J}=\sum_{i=1}^n\sum_{k=0}^{N-1}
\frac{h}{2}(u^{f,i}_{k+1})^TW_fu^{f,i}_{k+1}+
\frac{h}{2}(u^{m,i}_{k+1})^TW_mu^{m,i}_{k+1},\\
\text{such that } \braces{(R_N^i,x_N^i\Pi_N^i,\gamma_N^i)=(R_d,x_d^{A_i}(\theta),\Pi_d,\gamma_d)}_{i=1}^n,\\
\text{subject to discrete equations of motion
\refeqn{xkp}--\refeqn{Pikp},}
\end{gather*}
where $W_f,W_m\in\Re^{3\times 3}$ are symmetric positive definite
matrices.

Since we have neglected the gravitational interactions between the rigid bodies, the dynamics of the rigid bodies are only coupled through the terminal boundary conditions. If the parameter $\theta$ and the assignment $A$ are prescribed so that the terminal configurations of all rigid bodies are completely determined, then the optimal control problems for $n$ rigid bodies can be solved independently using the computational approach presented in Section \ref{sec:opt}. The formation cost is the summation of the resulting costs of each rigid body. Therefore, the optimal formation control problem for multiple rigid bodies consists in finding the optimal value of the parameter and the optimal assignment of rigid bodies among the $n!$ possible assignments.
This problem formulation is similar to the optimal formation reconfiguration problem presented in~\cite{JunMarObe.CDC06} except that we include the combinatorial assignment problem explicitly in this paper.

This requires combining the optimal control approach and the integer/combinatorial assignment optimization over the non-trivial dynamics of rigid bodies on $\SE^n$. We present a computational approach to this integrated optimal control problem.

\subsection{Optimal Control of Rigid Bodies on $\SE^n$}\label{subsec:optn}
We first solve the optimal formation control problem assuming that an assignment $A$ is pre-determined and fixed. We use a hierarchical optimal control approach~\cite{JunMarObe.CDC06}. Since the parameter $\theta$ completely defines the terminal configuration of all rigid bodies for the fixed assignment $A$, it also determines the corresponding cost by taking the sum of the cost of the optimal trajectories for each rigid body. Thus, the optimization problem is decomposed into an outer optimization problem to find the optimal value of $\theta$ that minimizes the total cost, and an inner optimization problem to find the optimal control forces and moments for the given value of $\theta$. This is a consequence of the fact that,
\begin{align}\label{eqn:hier0}
    \min_{u,\theta} \mathcal{J}(u,\theta)=\min_{\theta'} \braces{\min_u\braces{\mathcal{J}(u,\theta)|\theta=\theta'}}.
\end{align}
The inner optimization problem is solved by using the computational approach given in Section \ref{sec:opt}. The optimal value of $\theta$ is found by using a parameter optimization method with an explicit expression for the gradient.

\textit{Sensitivity of the cost: } Based on the solution of the optimal control problem of the $i$-th body, the sensitivity of the cost with respect to the parameter can be obtained as follows. Let $c^i\in\Re$ be the contribution of the $i$-th body to the performance index so that $\mathcal{J}=\sum_{i=1}^n c^i$.
\begin{align}
    c^i=\sum_{k=0}^{N-1}\frac{h}{2}(u^{f,i}_{k+1})^TW_fu^{f,i}_{k+1}+
\frac{h}{2}(u^{m,i}_{k+1})^TW_mu^{m,i}_{k+1}.\label{eqn:ci}
\end{align}
Suppose that the variation of the initial condition and the terminal boundary condition are given by $z_0^i,z_N^i\in\Re^{12}$. Using \refeqn{varopt}, the corresponding variation of the initial multiplier for an optimized solution is given by
\begin{align}
    \delta \lambda_0^i = (\Phi_N^{12,i})^{-1}(-\Phi_N^{11,i} z^i_0+z^i_N).\label{eqn:dellam0}
\end{align}
Substituting this into \refeqn{varopt}, we obtain the
variation of the multiplier as
\begin{align}\label{eqn:dellamk}
    \delta \lambda_k^i & = \Phi_k^{21,i}z_0^i +
    \Phi_k^{22,i}\delta\lambda_0^i,\nonumber\\
    & =\Phi_k^{21,i}z_0^i +
    \Phi_k^{22,i}(\Phi_N^{12,i})^{-1}(-\Phi_N^{11,i} z_0^i+z_N^i).
\end{align}
Since the control input is expressed in terms of the multiplier by the optimality condition \refeqn{ufkp},\refeqn{umkp}, the costs given in \refeqn{ci} are represented as follows.
\begin{align*}
    c^i & = \frac{h}{2} \sum_{k=0}^{N-1} (\lambda^{2,i}_k)^T W_f^{-1}(\lambda^{2,i}_k)
    +(\lambda^{4,i}_k)^T W_m^{-1}(\lambda^{4,i}_m),\\
    & = \frac{h}{2} \sum_{k=0}^{N-1} (\lambda^i_k)^T W \lambda^i_k,
\end{align*}
where $\lambda^i_k=[\lambda^{1,i}_k;\lambda^{2,i}_k;\lambda^{3,i}_k;\lambda^{4,i}_k]\in\Re^{12}$ and $W=\mathrm{diag}[0_{3\times 3}, W_f^{-1}, 0_{3\times 3}, W_m^{-1}]\in\Re^{12\times 12}$. Using \refeqn{dellamk}, the variation
of the cost is given by
\begin{align}
    \delta c^i & = h\sum_{k=0}^{N-1} (\lambda_{k}^i)^T W
        \delta \lambda_{k}^i,\nonumber\\
        & =\bracket{h\sum_{k=0}^{N-1} (\lambda_{k}^i)^T W
        (\Phi_k^{21,i}-\Phi_k^{22,i}(\Phi_N^{12,i})^{-1}\Phi_N^{11,i})}z_0^i\nonumber\\
        &\quad +\bracket{h\sum_{k=0}^{N-1} (\lambda_{k}^i)^T W
        \Phi_k^{22,i} (\Phi^{12,i}_N)^{-1}}z_N^i.\label{eqn:delci}
\end{align}
These represent the sensitivities of the optimal costs with respect
to the initial condition and the terminal boundary
condition; $\deriv{c^i}{z_0^i}, \deriv{c^i}{z_N^i}\in\Re^{1\times 12}$.

\textit{Computational approach:}
The sensitivity of the performance index with respect to the target space parameter is given by
\begin{align}
    \deriv{\mathcal{J}}{\theta}=\sum_{i=1}^n\deriv{c^i}{\theta} = \sum_{i=1}^n\deriv{c^i}{x_N^i} \deriv{x_d^{A_i}}{\theta},\label{eqn:delJ}
\end{align}
where $\deriv{c^i}{x_N^i}\in\Re^{1\times 3}$ is composed of the fourth to sixth elements of $\deriv{c^i}{z_N^i}$, i.e. $\deriv{c^i}{x_N^i}=\deriv{c^i}{z_N^i}B$, where $B=[0_{3\times 3},I_{3\times 3},0_{3\times 3},0_{3\times 3}]^T\in\Re^{12\times 3}$.

A quasi-Newton method can be applied to the outer optimization using this gradient computation. For example, the Broyden-Fletcher-Goldfarb-Shanno (BFGS) method to solve an unconstrained nonlinear optimization problem with an approximated Hessian is summarized as follows~\cite{Kel.BK95}.

\vspace*{0.1cm} {
\renewcommand{\theenumi}{\arabic{enumi}}
\renewcommand{\labelenumi}{\theenumi:}
\newcommand{\tab}{\hspace*{0.6cm}}
\hrule\vspace*{0.08cm}
\begin{enumerate}
\item Guess an initial parameter $\theta$.
\item Solve the $n$ optimal control problems.
\item Find $\deriv{\mathcal{J}}{\theta}$ using \refeqn{delJ}.
\item \textbf{while} $\norm{\deriv{\mathcal{J}}{\theta}} > \epsilon$.
\item \tab Find a line search direction; $D=-H^{-1}\deriv{\mathcal{J}}{\theta}$.
\item \tab Perform line search $\theta=\theta+\alpha D$.
\item \tab Update the Hessian $H$.
\item \tab Solve the $n$ optimal control problems.
\item \tab Find $\deriv{\mathcal{J}}{\theta}$ using \refeqn{delJ}.
\item \textbf{end while}
\end{enumerate}
\vspace*{0.08cm} \hrule} \vspace*{0.1cm} \noindent
Here $\epsilon$ denotes a stopping
criterion and $\alpha$ is a scaling factor, respectively. The Hessian can be initialized with $H=I_{l\times l}$ and updated during the BFGS iterations.

The major computational burden is in the third step and the fifth step. The computation time for the numerical iteration presented in Section \ref{sec:opt} can be substantially reduced if we have a good guess of the initial multiplier. At each iteration, we store the optimized initial multiplier for the terminal boundary conditions, and we use the accumulated data to initialize the initial multiplier at the next iteration. For example, \refeqn{dellam0} can be used to find an educated guess of the initial multiplier for neighboring boundary conditions. This reduces the computational burden as the iterations proceed, which will be shown by a numerical example in Section \ref{sec:ne}.

\subsection{Assignment Optimization Problem}\label{subsec:comb}
Now, we solve the optimal formation control problem assuming that the target parameter $\theta$ is determined and fixed. For the given value of $\theta$, the $n$ desired points $\braces{x_d^i}_{i=1}^n$, at which all rigid bodies are located at the terminal maneuver time, are completely defined. Thus, there are $n!$ possible combinatorial assignments.

Let $\braces{c^{ij}}$ be a $n\times n$ matrix, referred to as a cost matrix. Each element $c^{ij}$ represents the optimal cost of the $i$-th rigid body transferred to the $j$-th desired location. For an assignment $\braces{a_{ij}}$, the performance index is given by $\mathcal{J}=\sum_{i,j=1}^n c^{ij}a_{ij}$. The optimal assignment problem is given by
\begin{gather*}
    \min_{a_{ij}} \sum_{i,j=1}^n c^{ij} a_{ij},\\
    \begin{aligned}
    \text{Subject to   } &\sum_{j=1}^n a_{ij} =1 \quad \text{for } i\in\braces{1,2,\ldots,n},\\
                      &\sum_{i=1}^n a_{ij} =1 \quad \text{for } j\in\braces{1,2,\ldots,n},\\
                      &a_{ij} \in  \braces{0,1} \quad \text{for } i,j\in\braces{1,2,\ldots,n}.
    \end{aligned}
\end{gather*}
Since we assume that there is no interaction between rigid bodies, the cost matrix is independent of the assignment. For the given value of the target parameter $\theta$, we must solve at most $n^2$ optimal control problems to obtain the cost matrix. Once we have the complete cost matrix, the optimal assignment can be obtained by comparing costs for all possible assignments or by using the Hungarian method for large dimensional systems~\cite{Mur.BK85}.

It is often expensive to obtain the cost matrix. Each element of the cost matrix is a solution of the optimal control problem presented in Section \ref{sec:opt}. For the formation optimization problem, we need to find the cost matrix with varying values of the target parameter $\theta$. Thus, the objective of this subsection is to find the optimal assignment \textit{without solving} all $n^2$ optimal control problems. We start with an initial single spacecraft optimal trajectory computation, and use its optimal cost and the sensitivities given in \refeqn{delci} to populate the remaining entries of the cost matrix.

Suppose that we solve the optimal control problem of the first rigid body transferred to the first desired point to obtain $c^{11}$. Since this optimal solution is obtained by computing the linearized equations in \refeqn{varopt}. we can find the sensitivity of the cost with respect to the terminal boundary condition $\deriv{c^{11}}{x_d}$ by using \refeqn{delci}, without need of additional computational burden. Then, the optimal cost of transferring the first body to the other desired points, say $x_d$, is approximated as
\begin{align}
\hat c^{1}(x_d) = c^{11} + \deriv{c^{11}}{x_d} \Delta x_d
+ \frac{1}{2} (\Delta x_d)^T\frac{\partial^2 c^{11}}{\partial (x_d)^2}\Delta x_d,\label{eqn:c1approx}
\end{align}
where $\Delta x_d=x_d-x_d^1\in\Re^3$. The first order sensitivity is computed exactly, and the second order Hessian $\frac{\partial^2 c^{11}}{\partial x^1_d}$ is initially set to zero, and the approximation is improved as other optimal solutions become available. For example, if the solution of the optimal control problem of the first rigid body transferred to the second desired point is found, we obtain the exact value of $c^{12}$ and $\deriv{c^{12}}{x_d}$. This provides the following 1 and 3 dimensional constraints on the Hessian,
\begin{align*}
    \hat c^{1}(x_d^2)=c^{12},\quad\deriv{\hat c^{1}}{x_d}\bigg|_{x_d=x_d^2}=\deriv{c^{12}}{x_d}.
\end{align*}
Thus, the 6 elements of the Hessian can be approximated in either the minimum norm or least squares sense if additional optimal solutions involving the first rigid body are available.

This approach approximates the elements of the cost matrix along rows using the sensitivity of the cost with respect to terminal boundary conditions. A similar approximation along columns can be made by using the sensitivity of the cost with respect to the initial conditions. In the combinatorial optimization process, we utilize both approximations in order to avoid local minima. The advantage is that we use all of the sensitivity information available up to the current iteration in order to estimate the new cost matrix.

We construct a combinatorial assignment method using these approximations.
\renewcommand{\theenumi}{\roman{enumi}}
\begin{enumerate}
\item Guess an initial assignment and solve the corresponding
optimal control problems for this assignment.
\item Estimate the cost matrix using the linear approximations.
\item Find a new assignment using the estimated cost
matrix, and solve the corresponding optimal control problems for this assignment.
\item Find the best assignment using all of the solutions of the optimal
control problems obtained so far.
\item Construct a second-order approximation equation \refeqn{c1approx}
based on the best assignment, and estimate the cost matrix.
\item Repeat (iii)-(v) until the same approximation is repeated $M$ times in a row at (v).
\end{enumerate}
Numerical simulations show that setting $M=3$ is sufficient to find the global optimal assignment. At Steps (ii) and (v), we construct the approximation to the cost matrix using either the sensitivity of the cost with respect to the initial conditions or the sensitivity of the cost with respect to the terminal conditions. Several ways are summarized in Table \ref{tab:sensi}. Numerical simulations demonstrate that using both types of sensitivities has advantages and the last method is more efficient than the others in terms of finding the global optimal assignment with minimal computational effort.

\begin{table}
\renewcommand{\arraystretch}{1.3}
\caption{Methods to choose sensitivities}\label{tab:sensi}
\begin{tabular}{|c|p{7cm}|}\hline
Method & Sensitivity Selection\\\hline
Term.& Use terminal sensitivity always.\\
Init.& Use initial sensitivity always.\\
Rand.& Choose one of sensitivities randomly.\\
Rpt.& If the same assignment is repeated two times, switch to the different
sensitivity.\\
Alt.& Alternate sensitivities.\\
Comp.& For each element of the cost matrix, compare the
number of available solutions along the row direction, and the
number of available solutions along the column direction. Select one
of sensitivities that has more available
solutions. If both directions have the same number of available
solutions, choose randomly.\\\hline
\end{tabular}
\end{table}

\subsection{Computational Approach for Optimal Formation Control}
We have presented two optimization approaches; finding the optimal value of the space parameter for a given assignment, and finding the optimal assignment for a given value of the target parameter. We integrate both methods using a hierarchical optimization approach similar to \refeqn{hier0}.

The original optimization problem is stated as finding the optimal control forces, moments, target parameter, and assignment that minimizes the total cost.
\begin{align*}
    \min_{u,\theta, A}\mathcal{J} (u,\theta, A).
\end{align*}
Equivalently, this can be stated as finding the optimal assignment over the solutions for the optimal control inputs and the target parameters
\begin{align*}
    \min_{A'}\braces{\min_{u,\theta} \braces{\mathcal{J} (u,\theta, A)|A=A'}}.
\end{align*}
In the inner stage, we optimize the target parameter and the control inputs using the continuous optimization approaches presented in Section \ref{subsec:optn}, and in the outer stage, we find the optimal assignment using the combinatorial optimization approach presented in Section \ref{subsec:comb}. The optimization process is terminated when the iterations yield a solution that is optimal for both the continuous and combinatorial optimization stages.

\section{Numerical Example}\label{sec:ne}

We study a maneuver involving 5 identical rigid spacecraft under a central gravity
field. We assume that the mass of each spacecraft is negligible
compared to the mass of a central body, and we consider a fixed
frame attached to the central body as an inertial frame. The
resulting model is a Restricted Full Body Problem (RFBP)~\cite{CDC06.opt}.

Each spacecraft is modeled as a dumbbell, which consists of two equal
spheres and a massless rod. The gravitational potential is given by
\begin{align}
U=\sum_{i=1}^n -\frac{GMm}{2} \sum_{q=1}^2 \frac{1}{\norm{x^i+R^i\rho^i_q}},
\end{align}
where $G\in\Re$ is the gravitational constant, $M,m\in\Re$ are the
mass of the central body, and the mass of the dumbbell,
respectively. The vector $\rho^i_q\in\Re^3$ is the position of the
$q$th sphere from the mass center of the $i$-th dumbbell expressed in the
body fixed frame ($q\in \braces{1,2}$). The mass, length, and time
dimensions are normalized by the mass of the dumbbell, the radius of
a reference circular orbit, and its orbital period.

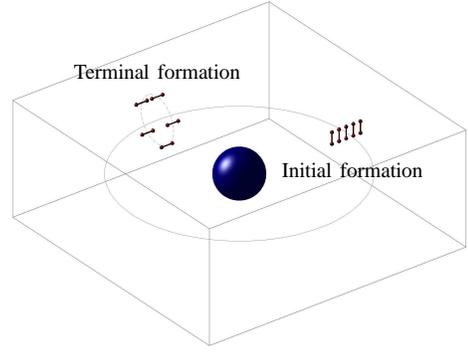
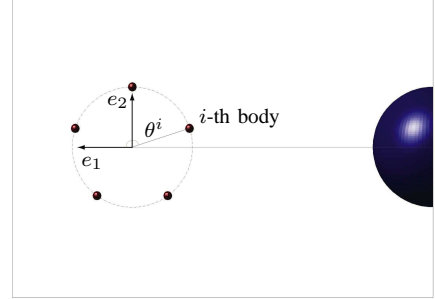
\begin{figure}[t]
\renewcommand{\xyWARMinclude}[1]{\includegraphics[width=0.70\columnwidth]{#1}}
\centerline{%
    \subfigure[The initial formation and the terminal formation]{\footnotesize\selectfont \begin{xy}%
    \xyWARMprocessEPS{formation}{eps}
    \xyMarkedImport{}
    \xyMarkedMathPoints{6,12}%
    \end{xy}\label{fig:forit}}}%
\centerline{%
\renewcommand{\xyWARMinclude}[1]{\includegraphics[width=0.65\columnwidth]{#1}}
    \subfigure[The terminal formation on a target circle]{\footnotesize\selectfont $$\begin{xy}%
    \xyWARMprocessEPS{formation2d}{eps}
    \xyMarkedImport{}
    \xyMarkedMathPoints{1-4}%
    \end{xy}$$\label{fig:theta}}%
} \caption{The initial formation and the desired terminal formation of 5 dumbbell
spacecraft on a circle}\label{fig:sco}
\end{figure}

The spacecraft are initially aligned along a radial direction as shown in \reffig{forit}. At the terminal time, we require that spacecraft are equally distributed on a target circle described by the location of its center $x_\circ\in\Re^3$, the radius $r_\circ\in\Re$, and the unit normal vector $n_\circ\in \S^2$. Let $\theta^i\in\S^1$ be the angle of the $i$-th spacecraft on the target circle from a given reference direction as shown in \reffig{theta}. We choose the target parameter as the angle of the first rigid body. The target $\mathcal{T}$ is given by
\begin{align*}
\mathcal{T}(\theta^1,A)=\braces{x_\circ+r_\circ\cos\theta^i e_1 + r_\circ \sin\theta^i e_2}_{i=1}^n,
\end{align*}
where $e_1=\frac{x_\circ}{\norm{x_\circ}}$, $e_2=e_1\times n_\circ$ are unit vectors in the target plane, and the angle $\theta^i$ is chosen to distribute the spacecraft uniformly on the circle
\begin{align*}
\theta^i(\theta^1,A)=\theta^1+\frac{2\pi}{5}(A_i-i).
\end{align*}
Since the target parameter $\theta^1$ determines the terminal position of the first spacecraft completely, we require that the first spacecraft be assigned to the first desired location, i.e. $A_1=1$. There remains $4!$ assignments for the other four spacecraft. Thus, the optimization parameters are the angle of the first spacecraft on the target circle, the $4!$ assignments for the remaining spacecraft, and the control inputs and moments.

\begin{figure}[t]
    \centerline{
    \includegraphics[width=0.90\columnwidth]{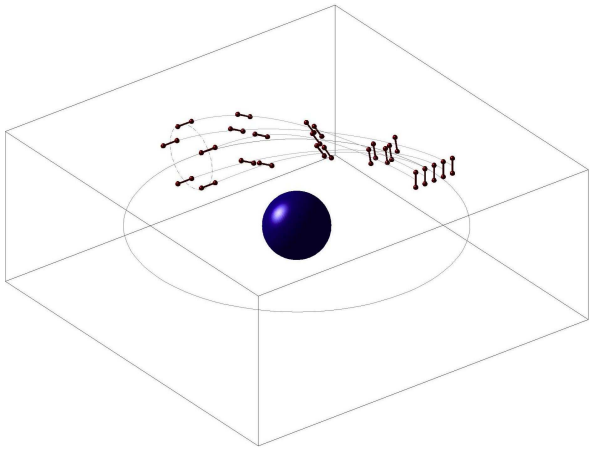}}
    \caption{Optimal spacecraft formation reconfiguration maneuver}\label{fig:fig3d}
\end{figure}

The iteration procedure for a particular numerical implementation of the optimization are shown as follows. The target parameter, assignment, cost, and computation time on an Intel Pentium M 1.73GHz processor using MATLAB are given for each iteration step.
\begin{enumerate}
\item[i)] The initial guess of the assignment is given by $A=\braces{(1,1),(2,4),(3,2),(4,3),(5,5)}$.
\item[ii.a)] For the given assignment, the optimal value of $\theta^1=2.4520$ is obtained in
$48.32$ minutes with cost $\mathcal{J}=8.6984$.
\item[ii.b)] For the given value of $\theta^1$, the optimal assignment of
$A=\braces{(1,1),(2,5),(3,2),(4,3),(5,4)}$ is obtained in $3.04$ minutes with cost
$\mathcal{J}=8.6905$.
\item[ii.c)] For the given $\theta^1$ and the given assignment, we check
$\deriv{\mathcal{J}}{\theta^1}=2.42\times 10^{-2}.$ Repeat iteration.
\item[iii.a)] For the given assignment, the optimal value of $\theta^1=2.5084$ is obtained in $12.69$ minutes with cost $\mathcal{J}=8.6898$.
\item[iii.b)] For the given value of $\theta^1$, the same optimal assignment of
$A=\braces{(1,1),(2,5),(3,2),(4,3),(5,4)}$ is obtained in $2.98$ minutes with cost
$\mathcal{J}=8.6898$.
\item[iii.c)] For the given $\theta^1$ and the corresponding optimal assignment, we check
$\deriv{\mathcal{J}}{\theta^1}=9.99\times 10^{-5}.$
\item[iv)] The optimization is terminated with optimal cost $\mathcal{J}=8.6898$ for $\theta^1=2.5084$ and $A=\braces{(1,1),(2,5),(3,2),(4,3),(5,4)}$
in total computation time $67.03$ minutes.
\end{enumerate}
The corresponding maneuvers for all the spacecraft are shown in \reffig{fig3d}.
The computation time for optimizing the target space parameter is reduced from $48.32$ minutes at Step (ii.a) to $12.69$ minutes at Step (iii.a). At each iteration, we use the optimization data accumulated in the previous iterations in order to initialize the initial multiplier for the optimal control problems. This reduces the computation time as the iterations proceed.

\begin{figure}
    \centerline{
    \subfigure[Histogram of total costs for $2400$ solutions with varying target parameters and assignments]{\includegraphics[width=0.64\columnwidth]{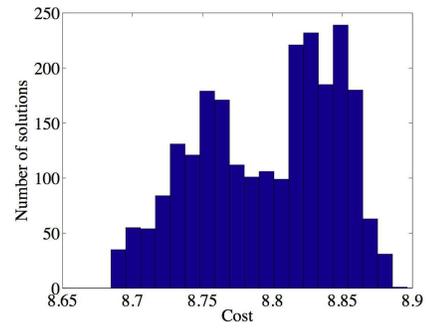}\label{fig:hist}}}
    \centerline{
    \subfigure[Histogram of total cost for $24$ optimized solutions with varying initial guesses of assignment]{\includegraphics[width=0.64\columnwidth]{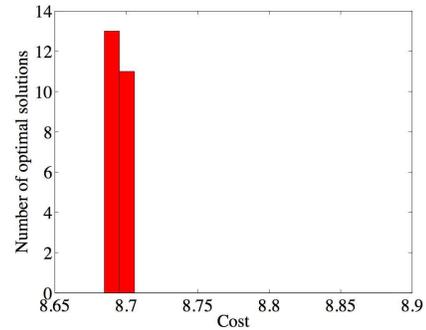}\label{fig:histopt}}}
    \caption{Distribution of the total costs before and after optimization}
\end{figure}

In order to estimate the distribution of the possible solutions, we uniformly discretize the interval $[0,2\pi)$ by 100 points for the target parameters and we find the total costs of $4!$ assignments for each value of the target parameter. The histogram for total costs of the corresponding $100\times 4!=2400$ solutions is shown in \reffig{hist}.

Numerical simulations show that the numerical optimized solution obtained depends on the initial guess of the assignment, and it is independent of the initial guesses for the target parameter and the initial Lagrange multiplier. We repeat the numerical optimization for all possible $4!$ initial guesses of the assignment. \reffig{histopt} shows the histogram of the optimized total costs for varying initial guesses for the assignment. Six initial assignments converged to the global optimal solution with $\mathcal{J}=8.6898$.

\section{Conclusions}
A combinatorial optimization method for spacecraft formation flight is presented. The objective is to transfer a group of identical spacecraft to a desired formation with minimum fuel expenditure. The assignment optimization over the discrete permutation group is explicitly integrated with the solutions of optimal control problems for combined orbital and rotational maneuvers of spacecraft, which are described by a Lie group variational integrator. The computational efficiency of the presented method is demonstrated by means of a numerical example.

\bibliography{opt}

\begin{thebibliography}{10}
\providecommand{\url}[1]{#1}
\csname url@rmstyle\endcsname
\providecommand{\newblock}{\relax}
\providecommand{\bibinfo}[2]{#2}
\providecommand\BIBentrySTDinterwordspacing{\spaceskip=0pt\relax}
\providecommand\BIBentryALTinterwordstretchfactor{4}
\providecommand\BIBentryALTinterwordspacing{\spaceskip=\fontdimen2\font plus
\BIBentryALTinterwordstretchfactor\fontdimen3\font minus
  \fontdimen4\font\relax}
\providecommand\BIBforeignlanguage[2]{{%
\expandafter\ifx\csname l@#1\endcsname\relax
\typeout{** WARNING: IEEEtran.bst: No hyphenation pattern has been}%
\typeout{** loaded for the language `#1'. Using the pattern for}%
\typeout{** the default language instead.}%
\else
\language=\csname l@#1\endcsname
\fi
#2}}

\bibitem{WanHad.JAS99}
P.~K.~C. Wang and F.~Y. Hadaegh, ``Minimum-fuel formation reconfiguration of
  multiple free-flying spacecraft,'' \emph{Journal of the Astronautical
  Sciences}, vol.~47, no. 1-2, pp. 77--102, 1999.

\bibitem{SavBulFra.CDC06}
K.~Savla, F.~Bullo, and E.~Frazzoli, ``On traveling salesperson problems for a
  double integrator,'' in \emph{Proceedings of the IEEE Conference on Decision
  and Control}, San Diego, California, Dec 2006, pp. 5305--5310.

\bibitem{GuiSch.JGCD06}
V.~M. Guibout and D.~J. Scheeres, ``Spacecraft formation dynamics and design,''
  \emph{Journal of Guidance, Control, and Dynamics}, vol.~29, no.~1, pp.
  121--133, 2006.

\bibitem{HosAtk.AIAA05}
A.~B. Hoskins and E.~M. Atkins, ``Spacecraft formation optimization with a
  multi-impulse design,'' in \emph{Proceedings of the AIAA Guidance,
  Navigation, and Control Conference and Exhibit}, San Francisco, California,
  Aug 2005, {A}{I}{A}{A} 2005-5835.

\bibitem{CMAME05}
T.~Lee, M.~Leok, and N.~H. McClamroch, ``{L}ie group variational integrators
  for the {F}ull {B}ody problem,'' \emph{{C}omputer {M}ethods in {A}pplied
  {M}echanics and {E}ngineering}, 2005, accepted.

\bibitem{CDC06.opt}
------, ``Optimal control of a rigid body using geometrically exact
  computations on {S}{E}(3),'' in \emph{Proceedings of the IEEE Conference on
  Decision and Control}, San Diego, California, Dec 2006, pp. 2170--2175.

\bibitem{MarWes.AN01}
J.~E. Marsden and M.~West, ``Discrete mechanics and variational integrators,''
  \emph{Acta Numerica}, vol.~10, pp. 357--514, 2001.

\bibitem{IserMun.AN00}
A.~Iserles, H.~Z. Munthe-Kaas, S.~P. N{\o}rsett, and A.~Zanna, ``Lie-group
  methods,'' \emph{Acta Numerica}, vol.~9, pp. 215--365, 2000.

\bibitem{CMDA06}
T.~Lee, M.~Leok, and N.~H. McClamroch, ``Lie group variational integrators for
  the full body problem in orbital mechanics,'' \emph{Celestial Mechanics and
  Dynamical Astronomy}, 2006, submitted.

\bibitem{Bry.BK75}
A.~E. Bryson and Y.-C. Ho, \emph{Applied Optimal Control}.\hskip 1em plus 0.5em
  minus 0.4em\relax Hemisphere Publishing Corporation, 1975.

\bibitem{Kel.BK95}
C.~T. Kelley, \emph{Iterative Methods for Linear and Nonlinear
  Equations}.\hskip 1em plus 0.5em minus 0.4em\relax SIAM, 1995.

\bibitem{ACC07.opt}
\BIBentryALTinterwordspacing
T.~Lee, M.~Leok, and N.~H. McClamroch, ``Optimal attitude control for a rigid
  body with symmetry,'' in \emph{Proceedings of the American Control
  Conference}, New York, July 2007, accepted. [Online]. Available:
  \url{http://arxiv.org/abs/math.OC/06009482}
\BIBentrySTDinterwordspacing

\bibitem{JunMarObe.CDC06}
O.~Junge, J.~E. Marsden, and S.~Ober-Bl{\"o}baum, ``Optimal reconfiguration of
  formation flying spacecraft: a decentralized approach,'' in \emph{Proceedings
  of the IEEE Conference on Decision and Control}, San Diego, California, Dec
  2006, pp. 5210--5215.

\bibitem{Mur.BK85}
K.~G. Murty, \emph{Linear and combinatorial programming}.\hskip 1em plus 0.5em
  minus 0.4em\relax Wiley, 1985.

\end{thebibliography}
\bibliographystyle{IEEEtran}

\end{document}